\documentclass[11pt]{article}

\pdfoutput=1
\usepackage[pdftex]{color}
\usepackage{amssymb}
\usepackage{amsthm}
\usepackage{amsmath}
\usepackage{latexsym}
\usepackage{amscd}
\usepackage{graphicx}
\usepackage[pdftex, colorlinks=true, citecolor=green]{hyperref}
\usepackage{lscape}
\usepackage{setspace}
\usepackage{multirow}
\usepackage{enumerate}

\newcommand{\ben}{\begin{equation}}     
\newcommand{\eeqn}{\end{equation}}
\newcommand{\bey}{\begin{eqnarray}}
\newcommand{\eey}{\end{eqnarray}}

\newtheorem{thm}{Theorem}[section]

\newtheorem{conj}[thm]{Conjecture}

\newtheorem{defn}[thm]{Definition}

\begin{document}

\noindent {\Large
\textbf{Symbolic template iterations of complex quadratic maps}
}
\\\\
 Anca R\v{a}dulescu$^{*,}\footnote{Assistant Professor, Department of Mathematics, State University of New York at New Paltz; New York, USA; Phone: (845) 257-3532; Email: radulesa@newpaltz.edu}$, Ariel Pignatelli$^{2}$
\\
\indent $^1$ Department of Mathematics, SUNY New Paltz, NY 12561
\\
\indent $^2$ Department of Engineering,  SUNY New Paltz, NY 12561

\vspace{7mm}
\begin{abstract}

\noindent The behavior of orbits for iterated logistic maps has been widely studied since the dawn of discrete dynamics as a research field, in particular in the context of the complex family $f \colon \mathbb{C} \to  \mathbb{C}$, parametrized as $f_c(z) = z^2 + c$, with $c \in \mathbb{C}$. However, little is is known about orbit behavior if the map changes along with the iterations. We investigate in which ways the traditional theory of Fatou-Julia may still apply in this case, illustrating how the iteration pattern (symbolic template) can affect the topology of the Julia and Mandelbrot sets.

\vspace{2mm}
\noindent We briefly discuss the potential of this extension towards a variety of applications in genetic and neural coding,
since it investigates how an occasional or a reoccurring error in a replication or learning algorithm may affect the dynamic outcome.
\end{abstract}

\section{Introduction}

\subsection{Discrete dynamics of the logistic family}

The family of logistic maps has been over the years one of the most studied examples in the theory of discrete dynamical systems. In the context of real interval maps, typically parametrized as $f_\mu \colon [0,1] \to [0,1], \; f_\mu(x) = \mu x(1-x)$, for $\mu \in [0,4]$, results from kneading theory classify the possible orbits, and provide a relationship between the combinatorics of critical orbits and the complexity of the corresponding map  (measured, for example, via its topological entropy). In the context of complex functions $f_c \colon \mathbb{C} \to \mathbb{C}, \; f_c(z) = z^2+c$, for $c \in \mathbb{C}$, results going all the way back to the original theory of Fatou and Julia describe thoroughly the behavior of the orbits in the dynamic complex plane, as well as phenomena in the parameter plane.

On the practical side, discrete iterations in general, and the dynamics of quadratic functions in particular, have been used to model natural processes. For example, the quadratic family has been used for more than a decade to model integrate and fire neurons~\cite{ermentrout1986parabolic,brunel2003firing}, and iterations of simple discrete  maps are the ideal candidate for modeling replication in genetic algorithms.

However, in nature, it is unlikely that systems evolve according to the same identical dynamics along time. Rather, one expects that occasional, or even periodic errors may be made in the iteration process, or even that the iteration scheme may change in time, according to the system's new needs as it is adapting. Therefore, a more realistic mathematical context to model such phenomena is to consider time-dependant iterations, in which the iterated map may change between steps, evolve in time, or appear (with variable frequency) in conjunction with other maps in the iteration. Based on this premise, some recent studies have looked at how the traditional theory extends when alternating two different real~\cite{radulescu2005connected,radulescu2008computing} or complex maps~\cite{danca2013graphical,danca2009alternated} in the iteration. In the latter references, the authors introduce the extended Julia sets corresponding to the alternated iteration of two complex quadratic maps. They show, analytically and numerically, that the dichotomy in the classical single map case (where the Julia sets are either totally connected, for some values of $c$, or totally disconnected, for the remaining values of $c$) does not hold in this case, and that Julia sets for alternated maps can also be disconnected without being totally disconnected. They furhter relate this extension to the Fatou-Julia theorem in the case of complex polynomials of degree greater than two, and show that alternated Julia sets exhibit graphical alternation.

We further extend this idea by looking at iteration of two different functions, $f_{c_0}$ and $f_{c_1}$, according to a more general binary symbolic sequence, in which the zero positions correspond to iterating the function $f_{c_0}$ and the one positions correspond to iterating the function $f_{c_1}$. This generalization is a more appropriate representation of replication or learning algorithms that appear in nature, with patterns that evolve in time, and which may involve occasional, random or periodic ``errors''. In this paper we investigate, primarily from a visual and numerical perspective, the questions that  arise for such complex iterations, and to what extent existing results are expected to hold in this more general setting. We are currently working on a theoretical frame to support analytical proofs for the type of results conjectured here.


\subsection{The traditional Julia and Mandelbrot sets}

The prisoner set of a complex map $f$ is defined as the set of all points in the dynamic plane, whose orbits are bounded. The escape set of a complex map is the set of all points whose orbits are unbounded. The Julia set of $f$ is definded as their common boundary $J(f)$. The filled Julia set is the union of prisoner points with their boundary $J(f)$.

For polynomial maps, it has been shown that the connectivity of a map's Julia set is tightly related to the structure of its critical orbits (i.e., the orbits of the map's critical points). Due to extensive work spanning almost one centry, from Julia~\cite{julia1918memoire} and Fatou~\cite{fatou1919equations} until recent developments~\cite{branner1992iteration,qiu2009proof}, we now have the following result, known as the Fatou-Julia Theorem: 

\begin{thm}
For a polynomial with at least one critical orbit unbounded, the Julia set is totally disconnected if and only if all the bounded critical
orbits are aperiodic.
\end{thm}

\vspace{2mm}
\noindent It was later shown~\cite{carleson1993complex,devaney2006criterion} that, as a consequence for the degree two case, Julia sets are either connected (if the orbit of the critical point 0 is bounded) or totally disconnected (if the orbit of 0 is unbounded).

If one had hoped similar results to hold for periodic iterations of more than one map, Danca et al.~\cite{danca2009alternated} showed that this dichotomy is broken even at the level of a simple alternation of two quadratic maps: aside from the two classical situations, the alternation also produced disconnected Julia sets which were not totally disconnected. While for periodic iterations one may still be able to build upon existing work and relate the structure of Julia sets with the critical alternating orbits, for nonperiodic iterations, however, the foundation is unlikely to remain applicable, and one may have to construct entirely new methods.

\subsection{Our extension: definitions and notations}

As with the traditional Julia set, we will be working with the complex quadratic family 

$$\{ f_c \colon \mathbb{C} \rightarrow \mathbb{C} \; \slash \; f_c(z) =z^2 + c,  \text{ with } c \in \mathbb{C} \}$$ 

\noindent However, for each iteration process, we will be using a pair of maps in this family, $f_{c_0}$ and $f_{c_1}$, as follows:

\begin{defn}
Fix $c_0, c_1 \in \mathbb{C}$, and a sequence ${\bf s} = (s_n)_{n \geq0} \in {\cal L}(\{0,1\})$ (which we will call the \textbf{symbolic template} of the iteration). For any $\xi_0 \in \mathbb{C}$, the \textbf{template orbit} $o_{\bf s}(\xi_0) = (\xi_n)_{n \geq 0}$ is the sequence constructed recursively, for every $n \geq 0$, as:

$$\xi_{n+1} = f_{c_{s_n}} (\xi_n)$$
\end{defn}

\vspace{5mm}
\noindent In other words, for each $n \geq 0$, we iterate $f_{c_0} (z) = z^2 + c_0$ if the corresponding entry $s_n=0$ in the symbolic template, and we iterate $f_{c_1}(z) = z^2 + c_1$ otherwise (i.e., if $s_n=1$). The Julia set for a symbolic template system is naturally defined as an extension of the traditional one:

\begin{defn}
Fix $c_0, c_1 \in \mathbb{C}$. Then we define, for any ${\bf s} \in {\cal L}(\{0,1\})$:

\vspace{2mm}
The template prisoner set: $P_{c_0,c_1}({\bf s}) = \{ \xi_0 \in \mathbb{C} \; \slash \; o_{\bf s}(\xi_0) \text{ is bounded } \}$

\vspace{2mm}
The template escape set: $E_{c_0,c_1}({\bf s}) = \{ \xi_0 \in \mathbb{C} \; \slash \; o_{\bf s}(\xi_0) \text{ is not bounded } \}$

\vspace{2mm}
The template Julia set: $J_{c_0,c_1}({\bf s}) = \partial P_{c_0,c_1}({\bf s})  = \partial E_{c_0,c_1}({\bf s}) $

\end{defn}

\vspace{3mm}
\noindent Since, in order to construct a template Julia set, we need both  a complex parameter pair $(c_0,c_1)$ and a symbolic template, we can consider phenomena in two different parameter spaces: in $\mathbb{C}^2$ (for a fixed template ${\bf s}$) and in the template space of binary sequences $l(\{ 0,1 \})$ (for a fixed pair $(c_0,c_1)$). We define two types of ``Mandelbrot'' sets, as follows:

\begin{defn}
Fix ${\bf s} \in {\cal L}(\{ 0,1 \})$ symbolic sequence. The corresponding \textbf{fixed-template Mandelbrot set} is defined as:

$${\cal M}_{\bf s} = \{ (c_0,c_1) \in \mathbb{C}^2 \; \slash \; o_{\bf s}(0) \text{ is bounded } \}$$
\end{defn}

\begin{defn}
Fix $(c_0,c_1) \in \mathbb{C}^2$. The corresponding \textbf{fixed-maps Mandelbrot set} is defined as:

$${\cal M}_{\bf c_0,c_1} = \{ {\bf s} \in {\cal L}(\{ 0,1 \}) \; \slash \; o_{\bf s}(0) \text{ is bounded } \}$$
\end{defn}

\noindent Properties of the fixed template Mandelbrot set have been previously studied for period two templates (i.e., alternating maps). In ~\cite{danca2013graphical}, the authors illustrate their corresponding extension of the Mandelbrot set as a subset of the quaternion filed, as well as some of its two-dimensional cross-sections. In the context of our extension, one may be interested to investigate a few directions, among which:

\begin{itemize}

\item Assessing the topological (e.g., connectedness) and fractal properties (e.g., Hausdorff dimention of the boundary) of fixed template Mandelbrot sets ${\cal M}_{\bf s}$, and how these vary with changing the template. 

\item Understanding topological and measure theoretical properties of fixed map Mandelbrot sets ${\cal M}_{\bf c_0,c_1}$, and how there depend on changing the parameters $c_0$ and $c_1$.

\item Studying whether / which aspects of the Fatou-Julia theory still hold for this extension. For example, it is clearly no longer true that an unbounded template orbit of zero implies a totally disconnected template Julia set. However, is it possible that the Julia set remain connected if the template orbit of 0 is bounded? In which case , do we still have a result such as the following:

\end{itemize}

\begin{conj}
For any template ${\bf s}$, the fixed template Mandelbrot set ${\cal M}_{\bf s}$ s equivalent with:

$${\cal M}_{\bf s} = \{ (c_0,c_1) \in \mathbb{C}^2 \; \slash \; J_{c_0,c_1}({\bf s}) \text{ is connected } \}$$

For any $(c_0,c_1) \in \mathbb{C}^2$, the fixed map Mandelbrot set ${\cal M}_{\bf c_0,c_1}$  is equivalent with:

$${\cal M}_{\bf c_0,c_1} = \{ {\bf s} \in {\cal L}(\{ 0,1 \}) \; \slash \; J_{c_0,c_1}({\bf s}) \text{ is connected } \}$$

\end{conj}

\vspace{2mm}

\noindent To the best of our knowledge, properties of the fixed-map Mandelbrot set have never been studied, and we believe this to be the first time the concept is being used.\\

\noindent In the following section, we investigate iterations of complex qudratic maps according to a variety of templates, and we discuss in particular some questions that arise as natural extensions of the traditional problems. The paper is organized as follws: In Section~\ref{periodic} we discuss properties of our Julia and Mandelbrot sets for fixed periodic templates. We focus in particular on coupling an arbitrary quadratic map $f_c(z)=z^2+c$ with the trivial map $f_0(z) = z^2$ (whose Julia set is the unit circle).  We observe how the connectivity patterns of the fixed template Julia and Mandelbrot sets change when the template itself is modified, to incorporate the two iterations in a different periodic mixture. In Section~\ref{nonperiodic}, we examine the same properties for nonperiodic templates.  In Section~\ref{propagate} we view the insertion of the second map as an ``error'' in the iteration process of the first map, study the effects of such an error propagating along the template on the structure of the corresponding Julia set. In Section~\ref{fixed_map}, we view the fixed map Mandelbrot set as a subset of $[0,1]^2$, and we define hybrid Mandelbrot sets. Finally, in Section~\ref{discussion}, we briefly discuss potential applications of this theory to the modeling of natural systems.


\section{Periodic template dynamics}
\label{periodic}

\subsection{Periodic template Julia sets}

\noindent In this section, we will consider $k$-periodic templates: ${\bf s} = [s_1,...,s_k]$, with $s_{k+j} = s_k$, for all $j \leq 0$. In previous work, Danca et al. have considered iterations of alternating quadratic maps~\cite{danca2013graphical,danca2009alternated}. In our context, this corresponds to considering the two possible symbolic templates of period two ${\bf s} = [01]$ and ${\bf s} = [10]$. 

In the references, the authors show that some basic properties and results are still inherited from the traditional single map case, but also that some new results emerge. For example, for any alternating orbit, the sequence of even and odd iterates are simultaneously bounded or unbounded, so that the Julia set of the alternated maps $f_{c_0}$ and $f_{c_1}$ and the Julia set of the quartic map $f_{c_1} \circ f_{c_0}$ have the same connectivity type for any $c_0, c_1 \in \mathbb{C}$. The authors further investigate numerically, graphically and analytically the relationship between the connectivity type of Julia sets and the boundedness of the critical orbits of alternating maps, verifying the first part of the Fatou-Julia theorem in this case. On the other hand, they also show that the dichotomy (connected vs, totally disconnected) characteristic to classical Julia sets no longer applies under the extension. Other questions, such as conditions for local connectivity of the 4-dimensional Mandelbrot set or of 2-dimensional Mandelbrot slices, remain open to further investigation.

\begin{figure}[h!]
\begin{center}
\includegraphics[width=0.9\textwidth]{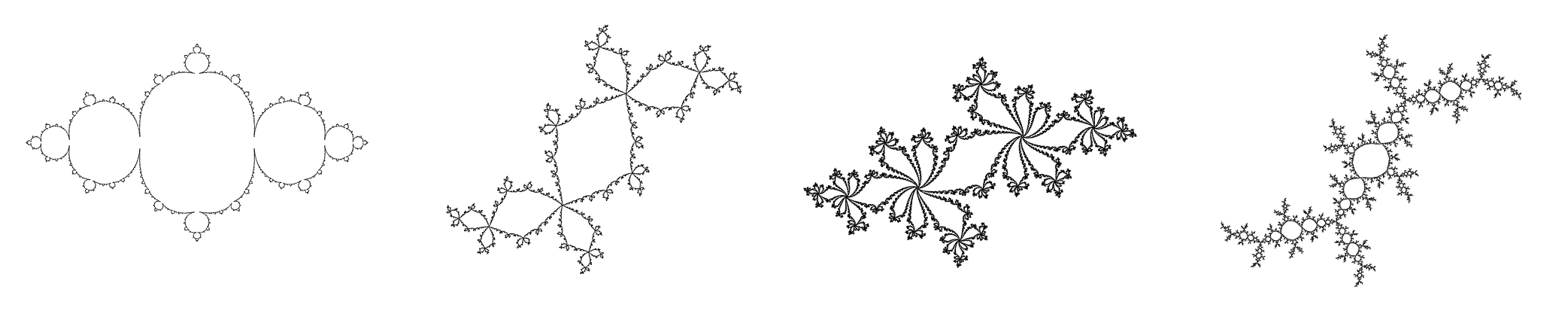}
\end{center}
\caption{\emph{{\bf Classical Julia sets} for four different values of the parameter $c$, for which the sets are connected: {\bf A.} $c=-0.75$; {\bf B.} $c=-0.117-0.76i$; {\bf C.} $c=-0.62-0.432i$; {\bf D.} $c=-0.117-0.856i$. We will be using these values in our subsequent illustrations, for combinations of two such functions iterated along a symbolic template.}}
\label{Julia_classic}
\end{figure}

Similar properties are likely to extend to periodic template iterations for periods $k>2$. Notice that a $j$-shift $\sigma$ on the periodic template $\sigma({\bf s}) = \sigma[s_1,...s_k] \to [s_{1+j},...s_{k+j}] = {\bf s}^j$ translates as a polynomial transformation of degree $j$ on the Julia set: 
$$J({\bf s}) = f_{c_{s_{j-1}}} \circ \hdots \circ f_{c_{s_1}} (J({\bf s}^j))$$

\noindent One may look at classes defined by shifts on periodic templates, so that Julia sets within the same class are all related by such transformations. The connectedness properties of a Julia set may be preserved or not by the transformation, depending on the parameters $c_0$ and $c_1$, and on the critical orbit properties produced by their combination.

For templates of period 3 for example, there are 3 such classes (or types): 

\begin{enumerate}[(1.)]
\setcounter{enumi}{-1}

\item All zero entries: ${\bf s} = [000]$

\item A single one entry: ${\bf s} = [001]$, ${\bf s} = [010]$, ${\bf s} = [100]$

\item Two one entries: ${\bf s} = [011]$, ${\bf s} = [101]$, ${\bf s} = [110]$

\item All one entries: ${\bf s} = [111]$

\end{enumerate}

While higher period template Julia sets show similar symmetry properties as observed in the alternating case~\cite{danca2013graphical}, new questions add to the discussion for higher periods. These refer to understanding how the structure of the periodic template block affects the resulting Julia set, \emph{in combination} with the parameter pair $(c_0,c_1)$. Connectivity of the Julia set is now influenced by three factors: (i) the maps $c_0$ and $c_1$, (ii) the balance of how often one map is iterated versus the other (i.e., the number of 1s vesus 0s in the template block) and (iii) the location of these 1s and 0s along the block. \\

\begin{figure}[h!]
\begin{center}
\includegraphics[width=0.9\textwidth]{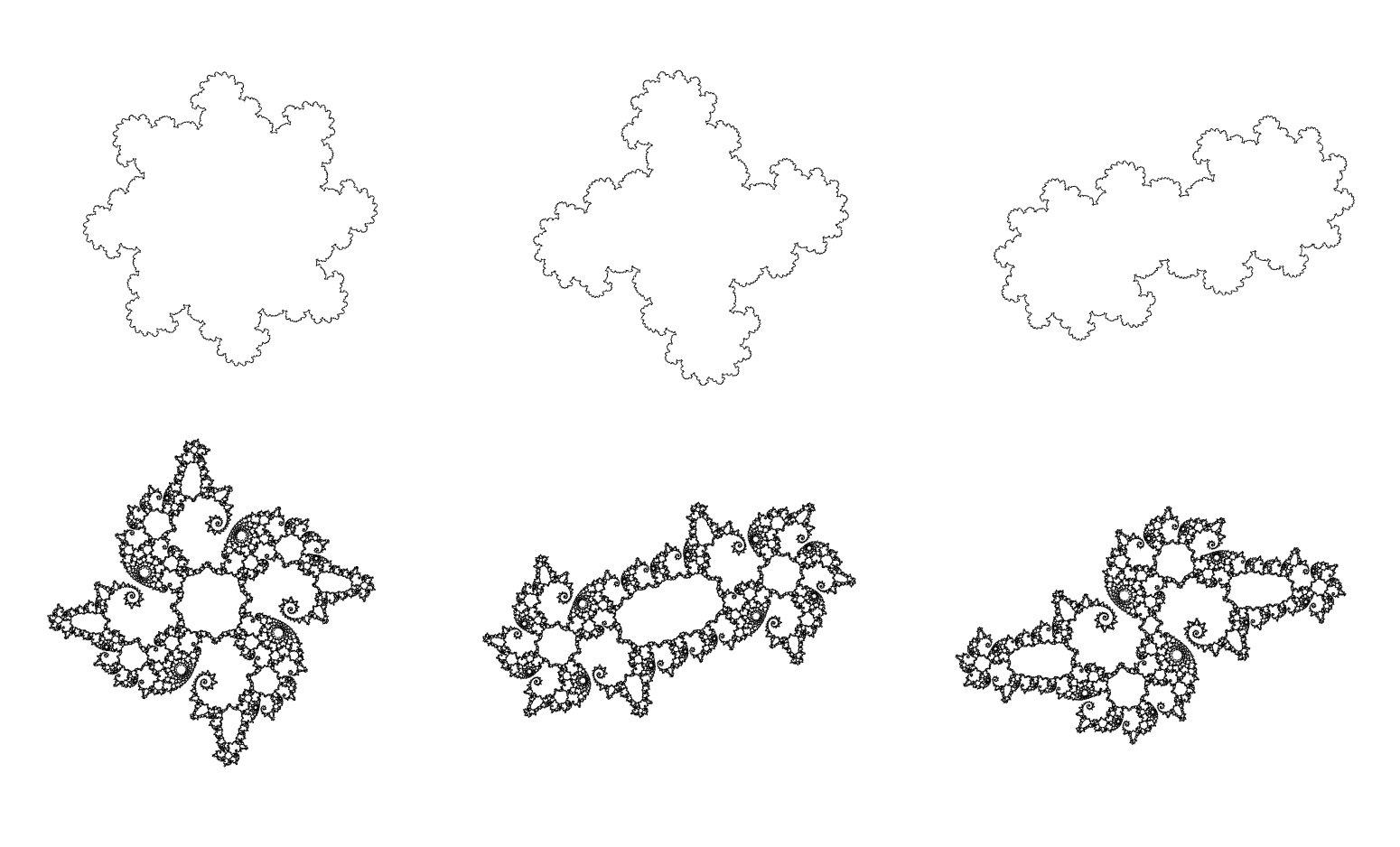}
\end{center}
\caption{\emph{{\bf Period 3 template Julia sets for $c_0 = 0$ and $c_1 = -0.62-0.432i$},  for templates of type 1 (top row) and type 2 (bottom row). {\bf A.} ${\bf s} = [011]$; {\bf B.} ${\bf s} = [101]$; {\bf C.} ${\bf s} = [110]$; {\bf D.} ${\bf s} = [001]$; {\bf E.} ${\bf s} = [010]$; {\bf F.} ${\bf s} = [100]$.}}
\label{Julia_c3}
\end{figure}

\begin{figure}[h!]
\begin{center}
\includegraphics[width=.8\textwidth]{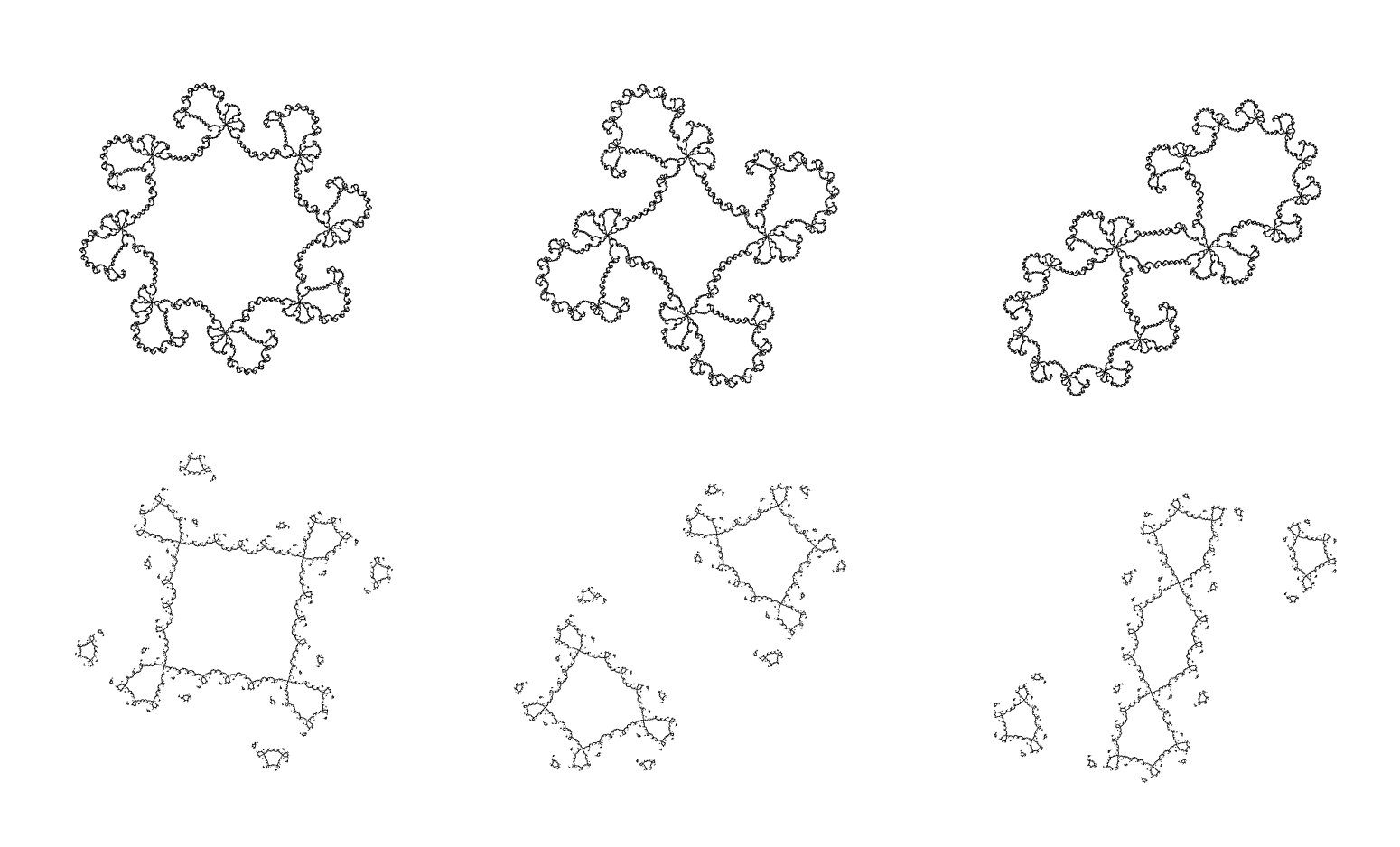}
\end{center}
\caption{\emph{{\bf Period 3 template Julia sets for $c_0 = 0$ and $c_1 = -0.117-0.76i$}, for templates of type 1 (top row) and type 2 (bottom row). {\bf A.} ${\bf s} = [011]$; {\bf B.} ${\bf s} = [101]$; {\bf C.} ${\bf s} = [110]$; {\bf D.} ${\bf s} = [001]$; {\bf E.} ${\bf s} = [010]$; {\bf F.} ${\bf s} = [100]$.}}
\label{Julia_c2}
\end{figure}

\begin{figure}[h!]
\begin{center}
\includegraphics[width=.9\textwidth]{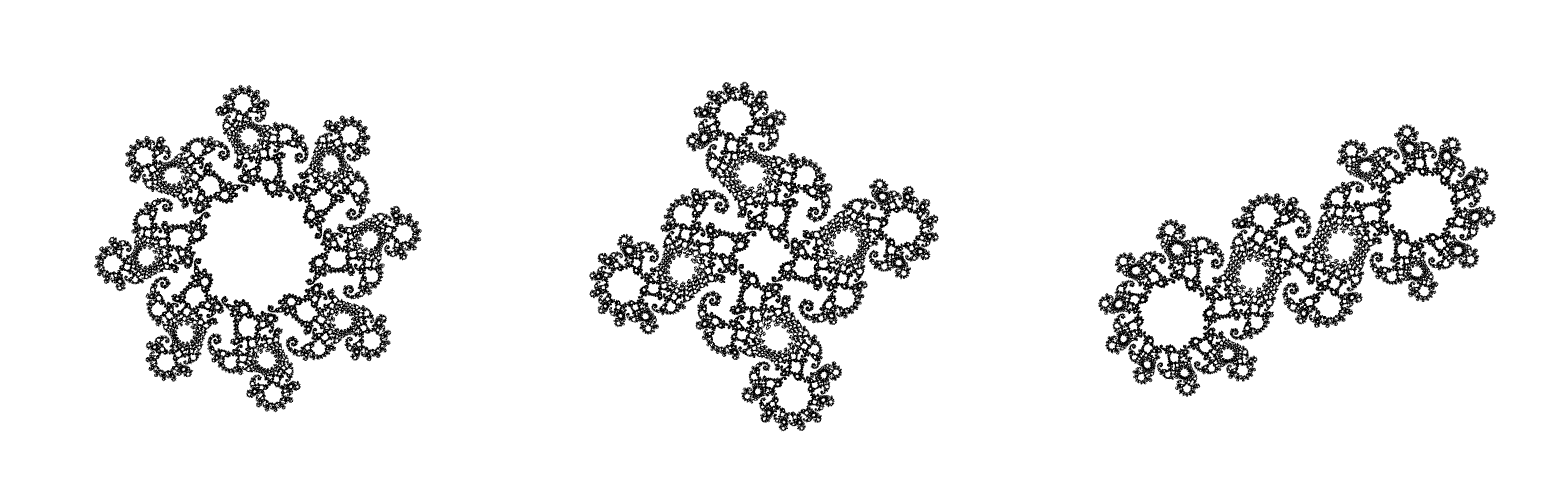}
\end{center}
\caption{\emph{{\bf Period 3 template Julia sets for $c_0 = 0$ and $c_1 = -0.5622-0.62i$}, for templates of type 1: {\bf A.} ${\bf s} = [001]$; {\bf B.} ${\bf s} = [010]$; {\bf C.} ${\bf s} = [100]$. The Julia sets for the same parameter values and templates of type 2, i.e. ${\bf s} = [011]$, ${\bf s} = [101]$, ${\bf s} = [110]$, are totally disconnected (dust).}}
\label{Julia_c4}
\end{figure}

\noindent Figures ~\ref{Julia_c3}, ~\ref{Julia_c2} and ~\ref{Julia_c4} show, for example, template Julia sets for 3-periodic templates, for combinations of the map $c_0 = 0$ with the maps $c_1 = -0.117-0.76i$, $c_1 = -0.62-0.432i$ $c_1 = -0.5622-0.62i$ and respectively $c_1 = -1-0.55i$. For the first two of these $(c_0,c_1)$ pairs, the classical Julia sets are shown in Figure~\ref{Julia_classic}; the other two Julia sets are totally disconnected.

\begin{figure}[h!]
\begin{center}
\includegraphics[width=0.9\textwidth]{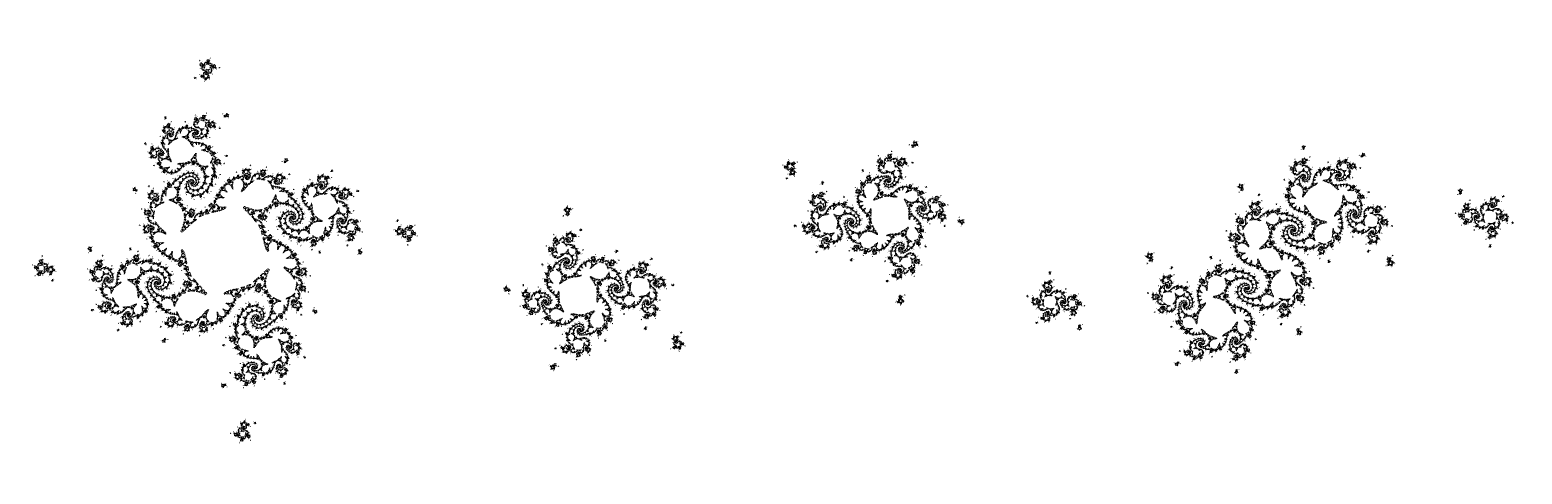}
\end{center}
\caption{\emph{{\bf Period 3 template Julia sets for $c_0 = 0$ and $c_1 =  -1-0.55i$}, for templates of type 2: {\bf A.} ${\bf s} = [011]$; {\bf B.} ${\bf s} = [101]$; {\bf C.} ${\bf s} = [110]$. The Julia sets for the same parameter values and templates of type 1, i.e. ${\bf s} = [001]$, ${\bf s} = [010]$, ${\bf s} = [100]$, are totally disconnected (dust).}}
\label{Julia_c15}
\end{figure}

The template type seems to affect connectivity of the template Julia set as significantly as the parameters $(c_0,c_1)$. Clearly, it is not surprising that, for fixed $c_0$ and $c_1$, using a different template type (e.g., 1 versus 2 ones in the 3-periodic case) will affect the connectivity of the Julia set. However, the results may be rather counterintuitive (see Figure~\ref{Julia_c15}). In addition to the template type, changing the position of the 1s along the template may or may not affect connectivity, as shown in Figures~\ref{Julia_c3} and~\ref{Julia_c2}  for the three different template blocks of type 2. In particular, the folding produced by applying a shift to the template may break or merge connectivity loci, depending on the template and on the complex parameters $c_0$ and $c_1$ used. 

For $c_1 = -0.62-0.432i$ (Figure~\ref{Julia_c3}), $c_1 = -0.117-0.76i$ (Figure~\ref{Julia_c2}) and $c_1 = -0.5622-0.62i$ (Figure~\ref{Julia_c4}) in combination with $c_0=0$, the Julia set becomes more disconnected as a result of increasing the contribution of $c_1$ versus $c_0$. Indeed, in the first case, the template Julia set, connected for all cases of type one template (a single one and two zeros) remains connected (although a lot more complex) for all corresponding type 2 templates (two ones and one zero). In the second case, while all type 1 templates render connected Julia sets, these become disconnected for all cases of type 2 template. In the third case, the Julia set for type 1 templates, already disconnected, turns to dust for all type 2 templates. This is somewhat expectable, given the shape of the Julia set for each of these maps considered in isolation (the unit circle for $c=0$, the Julia set shown in Figure~\ref{Julia_classic}c for $c=-0.62-0.432i$, in Figure~\ref{Julia_classic}b, for $c=-0.117-0.76i$, and respectively a totally disconnected Julia set for $c=-0.5622-0.62i$). However, for $c_1 =  -1-0.55i$ (Figure~\ref{Julia_c15}), the opposite happens: the template Julia sets are disconnected, but not totally disconnected for type 2 templates, and become dust for type 1 templates, where the more substantial contribution of the map $f_{c_0}(z) = z^2$ would suggest otherwise.

\subsection{Periodic template Mandelbrot sets}

We next observe the structure of template Mandelbrot sets for periodic templates. We are interested in particular to illustrate how the fractal structure of the set (represented by its Hausdorff dimension) and its connectedness properties change when changing the template density or regularity.

Since, as in the original reference, we are using two iterated maps in combination, the template Mandelbort set can be seen as a 2-dimensional complex object (the locus of $(c_0,c_1) \in \mathbb{C}^2$  for which the template orbit of zero is bounded). One can visualize the connectedness of the MAndelbrot set for a fixed template by looking at a latice of 2-dimensional slices, each representing the behavior with respect to $c_1$ for a fixed $c_0$ (see in Figure~\ref{Mand_all_011} in Appendix 2).

\begin{figure}[h!]
\begin{center}
\includegraphics[width=.8\textwidth]{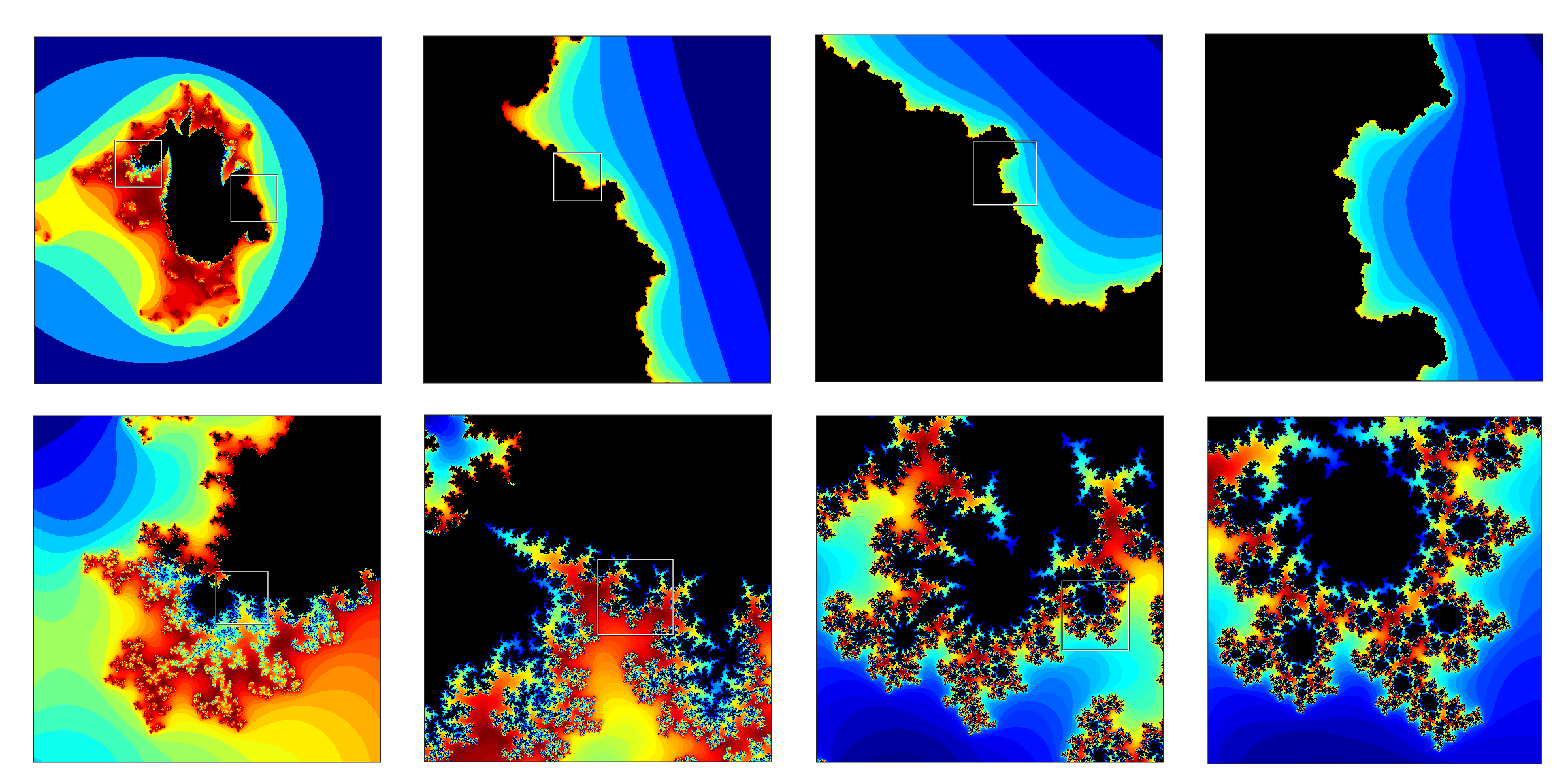}
\end{center}
\caption{\emph{{\bf Template Mandelbrot slice for the periodic template} [011], for fixed $c_0=-0.2+0.6i$, with two progressively zoomed in windows: one with low, the other with high Hausdorff dimension along the boundary. Three consecutive zoom-ins of the first instance are represented on top, and four zoom-ins of the second situation are represented on the bottom. The pattern persists at higher and higher levels, suggesting an alternation of high and low dimension intervals along the boundary of the set.}}
\label{Mand_zoom_011}
\end{figure}

\begin{figure}[h!]
\begin{center}
\includegraphics[width=.8\textwidth]{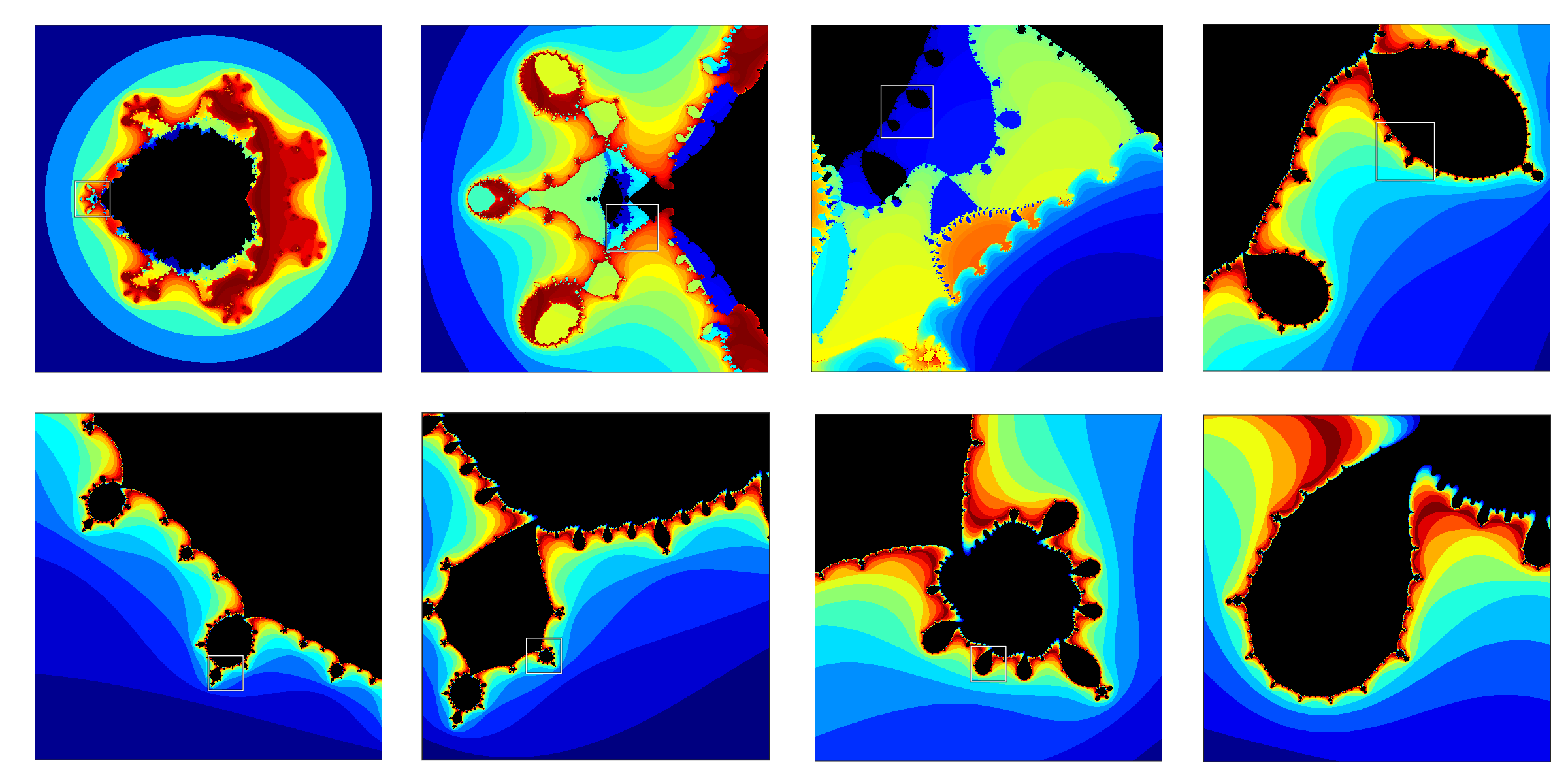}
\end{center}
\caption{\emph{ Template Mandelbrot set for a random template, for fixed $c_0=0$. Each panel shows at higher resolution a zoomed in window, marked on the preceding panel. The stucture persists at higher and higher levels, suggesting the preservation of fractal features in the Mandelbrot set, even for random templates.}}
\label{mand}
\end{figure}

When we tracked the fractality of the boundary of template Mandelbrot sets, we found consistent trends, which we describe visually in Figure~\ref{Mand_zoom_011}, for periodic templates, and in Figure~\ref{mand}, for nonperiodic templates, using consecutive zooms into various regions along the boundary.

we show consecutive zooms into the boundary structure, at deeper and deeper levels).

\section{Nonperiodic templates}
\label{nonperiodic}

The case of periodic templates represents only a first extension, with basic properties expected to replicate quite naturally the theory for iterations or alternated maps. We next consider the more general case, of binary templates which are not necessarily periodic. In this case, the existing methods cannot be applied directly. To begin with, in this more general case, there is no extension for the classic concept of critical orbits, hence no corresponding possibility of encoding the behavior of the Julia set in the behavior of these orbits. As a consequence, results such as an extension of the Fatou-Julia theorem, or local connectivity of fixed-template Mandelbrot sets may be harder to determine, and require new methods. In this paper we only describe the setup, draw comparisons with the classical case and lay out a few of the interesting questions and conjectures that arise in this new framework. A more rigorous analysis, discussed briefly in the last section, is the subject of our current research.

\subsection{Julia sets for nonperiodic templates}
\label{nonperiodic_Julia}

\begin{figure}[h!]
\begin{center}
\includegraphics[width=\textwidth]{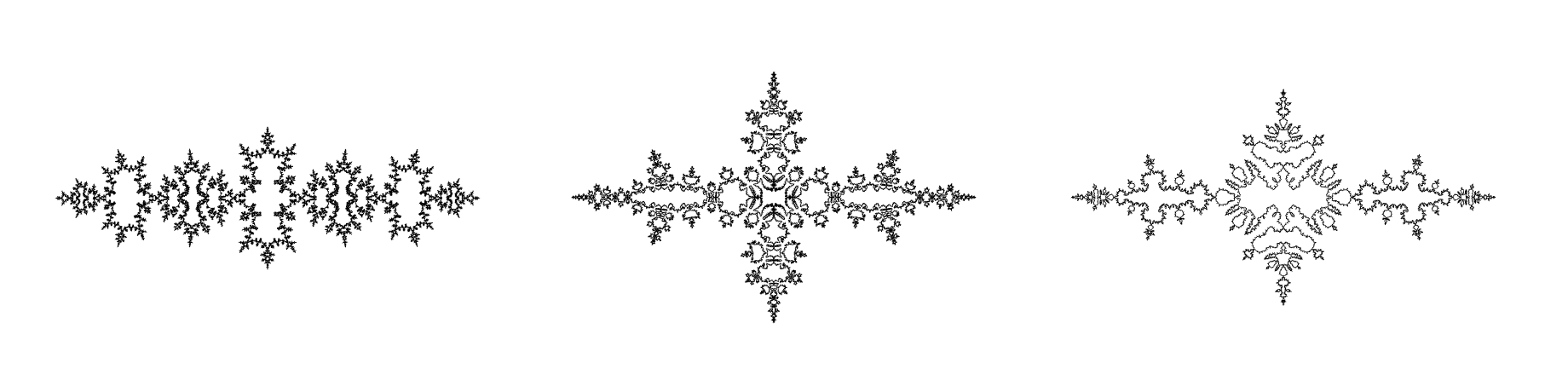}
\end{center}
\caption{\emph{{\bf Template Julia sets} for one fixed nonperiodic template, and different pairs $(c_0,c_1)$: {\bf A.} $c_1=-0.75$ and $c_0=-1.2$; {\bf B.} $c_1=-1.2$ and $c_0=-0.75$; {\bf C.} $c_1=-1.2$ and $c_0=-0.8$. For the simulation, the template creating the Julia set was truncated to 200 iterations.}}
\label{Julia_random_combined}
\end{figure}

In Figure~\ref{Julia_random_combined}, one nonperiodic template was used to create the Julia set for three different combinations of parameters. Perhaps the first thing one notices is that the Julia sets for random templates still exhibit complex (fractal) structure, symmetry and alternations. However, we need to recall that, in all numerically generated figures, one can only use truncated representations of infinite templates, retaining only a specific number of iterations when representing the Julia set (for example, in all our figures we used 200 iterations). A natural question concerns the variability with the number of iterations (i.e., length of the truncated template) of the resulting approximation for the Julia set.

\begin{figure}[h!]
\begin{center}
\includegraphics[width=.8\textwidth]{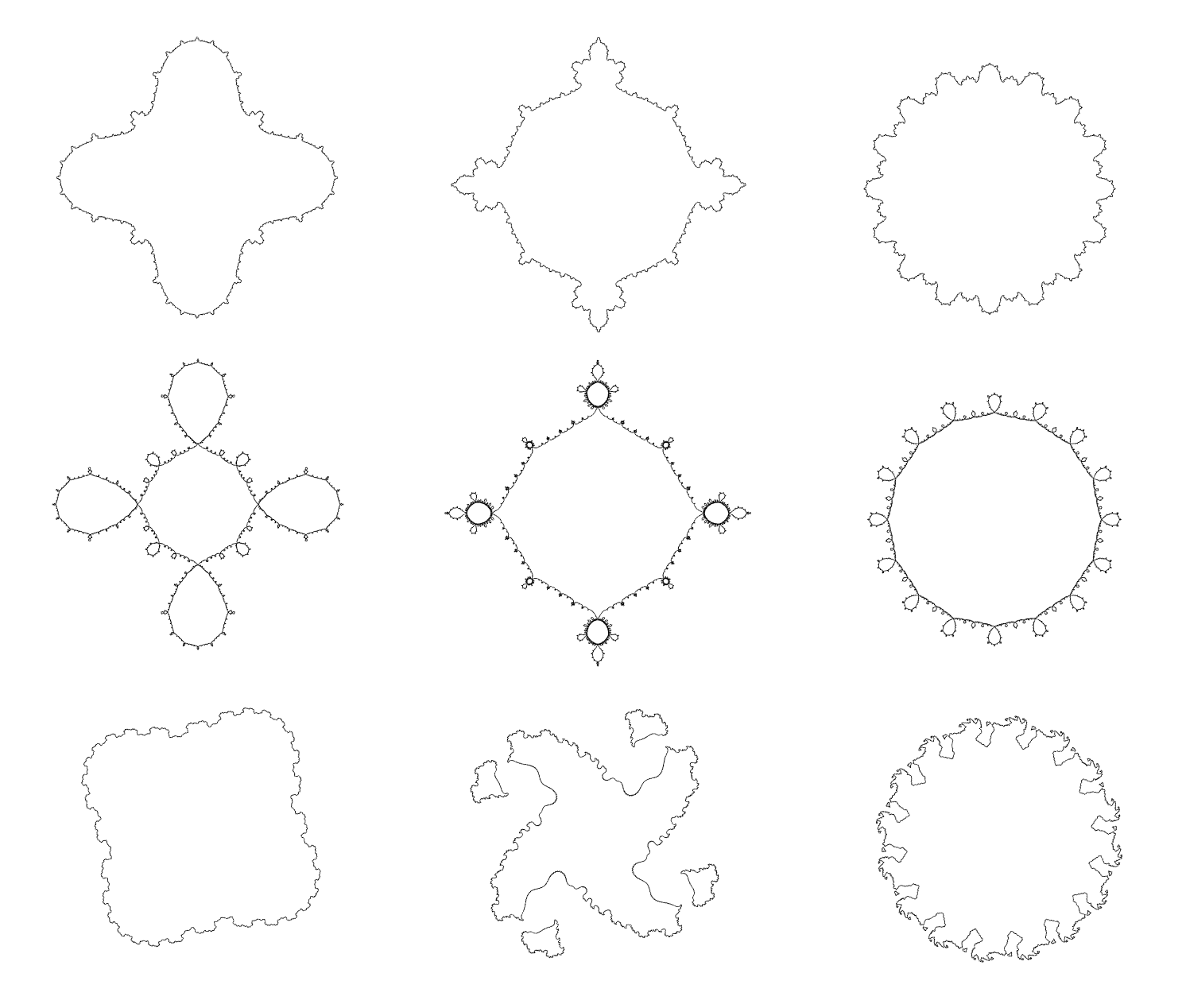}
\end{center}
\caption{\emph{{\bf Template Julia sets} for three different randomly generated templates (from left to right), with parameters $c_0=0$ and, from top to bottom: $c_1=-0.75$ (top), $c_1 = -0.8$ (middle)  and $c_1=0.375+0.333i$ (bottom). For the simulation, the template creating the Julia set was truncated to 200 iterations.}}
\label{Julia_random_c1_c8_c13}
\end{figure}

\begin{figure}[h!]
\begin{center}
\includegraphics[width=.8\textwidth]{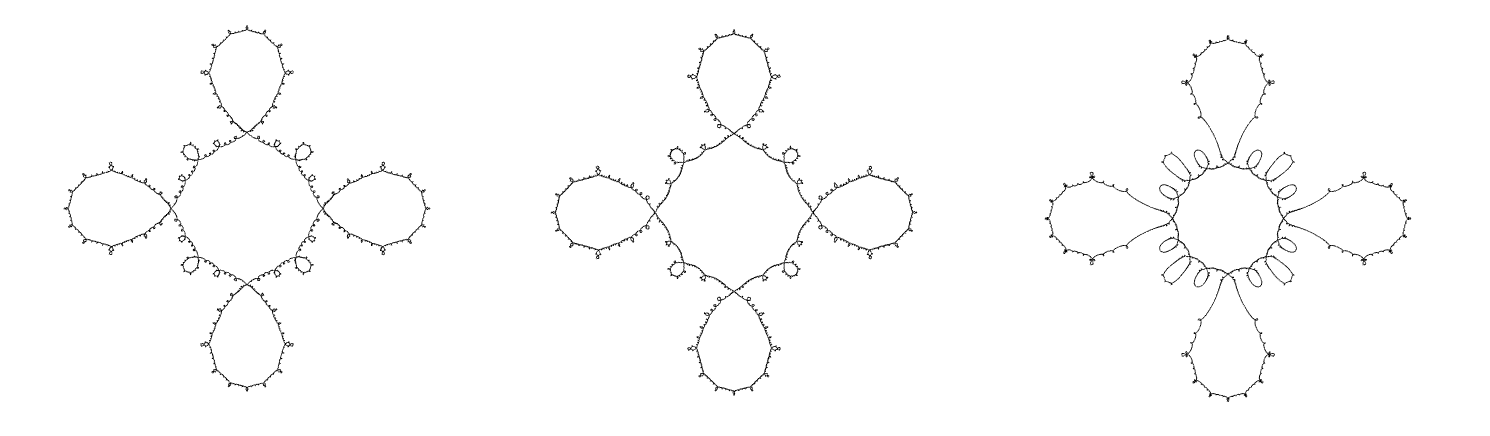}
\end{center}
\caption{\emph{{\bf Template Julia sets} for $c_1 = -0.8$ and $c_0=0$, for three random truncated templates of length 200 with the same root of length $l=10$. The first template (left) is the same as Template 1 in Figure~\ref{Julia_random_c1_c8_c13}, hence the corresponding Julia set is identical with that in the middle left panel in the referenced figure.}}
\label{Julia_common_root}
\end{figure}

In Figure~\ref{Julia_random_c1_c8_c13}, three different nonperiodic (randomly generated) templates were used to create the corresponding Julia sets, for the now familiar combination of parameter values $c_1 = -0.62-0.432i$ and $c_0=0$. Notice that different templates introduce differences in the Julia set as significant as changes in the parameters $(c_0,c_1)$. We then fixed Template 1 from Figure~\ref{Julia_random_c1_c8_c13} and $c_0=0$, $c_1=-0.8$ (that is, we started with the Julia set in the middle left panel). Using the same $(c_0,c_1)$ pair, we computed the Julia sets for a collection of templates of length 200 with identical entries up to the 10th position. The variations among the Julia sets we obtained this way are clearly much smaller, as illustrated in Figure~\ref{Julia_common_root}. To better understand this dependence, we consider the following:

\begin{defn}
For any $n \in \mathbb{N}$, we call the $k$-root of a template ${\bf s} = (s_k)$ the finite sequence $\overline{\bf s}^k = s_1, \hdots ,s_k$.
\end{defn}

\begin{defn}
We say that two templates have the same $k$-root if they agree up to their $k$-th position.
\end{defn}

\noindent Based on our numerical simulations (see Figure~\ref{Julia_common_root}), we propose the following conjectures:

\begin{conj}
Fix a parameter pair $(c_0,c_1)$, and a template ${\bf s}$.  For any $\delta > 0$, there exists $n \in \mathbb{N}$ such that, if ${\bf s}_1$ and ${\bf s}_2$ are two templates with a common $n$-root, then 
$$d(J_{c_0,c_1}({\bf s}_1),J_{c_0,c_1}({\bf s}_2)) < \delta$$
\noindent where $d$ represents the Hausdorff distance between two sets.
\end{conj}

\begin{conj}
Fix a parameter pair $(c_0,c_1)$. For any $\delta > 0$, there exists $N \in \mathbb{N}$ such that:
$$d(J_{c_0,c_1}(\overline{\bf s}^n) - J_{c_0,c_1}({\bf s})) < \delta \text{ , for } n \geq N$$
\end{conj}

\vspace{2mm}
\noindent In other words, the truncated Julia sets converge to the true Julia set in Hausdorff distance, as the root of the template is getting larger.\\

\subsection{Propagating a perturbation}
\label{propagate}

One way to interpret the combination of two maps in the iteration scheme is to think of $f_{c_1}$ as the desired map, and of $f_{c_0}$ as an ``error,'' or ``perturbation'' in the desired iteration. In this context, an all 1 template corresponds to a ``perfect'' replication  process, in which a specific map $f_{c_1}$ is iterated identically for a (large or infinite) number of times, and any number of zero entries in the template correspond to as many intrusions of the erroneous map $f_{c_0}$ in this iteration process. This intrusion can be periodic (as discussed in Section~\ref{periodic}), or may occur at a random sequence of steps (as discussed in Section~\ref{nonperiodic_Julia}). In this section, we illustrate the changes in the Julia set produced by propagating a single error along a perfect template. 

In Appendix 1, we illustrate the effects on the Julia set of a single propagting error, for a few combinations of maps. In Figures~\ref{error_anim_c3} and~\ref{error_anim_c7}, the desired map $f_{c_1}$ is combined with the trivial map $f_0(z) = z^2$, while in Figure~\ref{error_anim_c7_c4}, the error is a different (and unrelated) function, with complex parameter $c_0 \neq 0$. In Figure~\ref{error_anim_c7_perturb}, the error is a small perturbation of the original function, i.e. $c_0 = c_1 + \epsilon$, with a small $\varepsilon \in \mathbb{C}$.

In general, we notice that the impact of the same error on the Julia set is more substantial when the error occurs earlier in the iteration process. In fact, the Julia sets start to be reminiscent of the correct Julia set less than 50 iterations through the template length. Along the way, however, the connectivity and symmetry of the erroneous Julia changes at each step, sometimes in unexpected and counter-intuitive ways. Notice, for example, that the small perturbation to $c_1$ (in Figure~\ref{error_anim_c7_perturb}) has at each step an effect on the Julia set which is comparable with that of a much larger perturbation (in Figures~\ref{error_anim_c7} and~\ref{error_anim_c7_c4} ). This promotes the possibility that, in such a replication process, the timing of the error is equally or even more important than the size of the error.\\


\subsection{Mandelbrot sets for nonperiodic templates}
\label{nonperiodic_Mandelbrot}

Due to the lack of regularity in the occurences of $c_0$ and $c_1$ in the template, the fractal dimension of the template Mandelbrot set seems to be lower than that for periodic templates. We conjecture, however, that the Hausdorff dimension is more uniform along the boundary of a set whose template exhibits less regularity -- in contrast with the situation for regular templates. There, the fractal behavior was much more variable along the boundary of the corresponding Mandelbrot set, with portions with very high Hausdorff dimension, and portion with low dimension. 

\section{Fixed-map Mandelbrod sets}
\label{fixed_map}

\begin{figure}[h!]
\begin{center}
\includegraphics[width=0.8\textwidth]{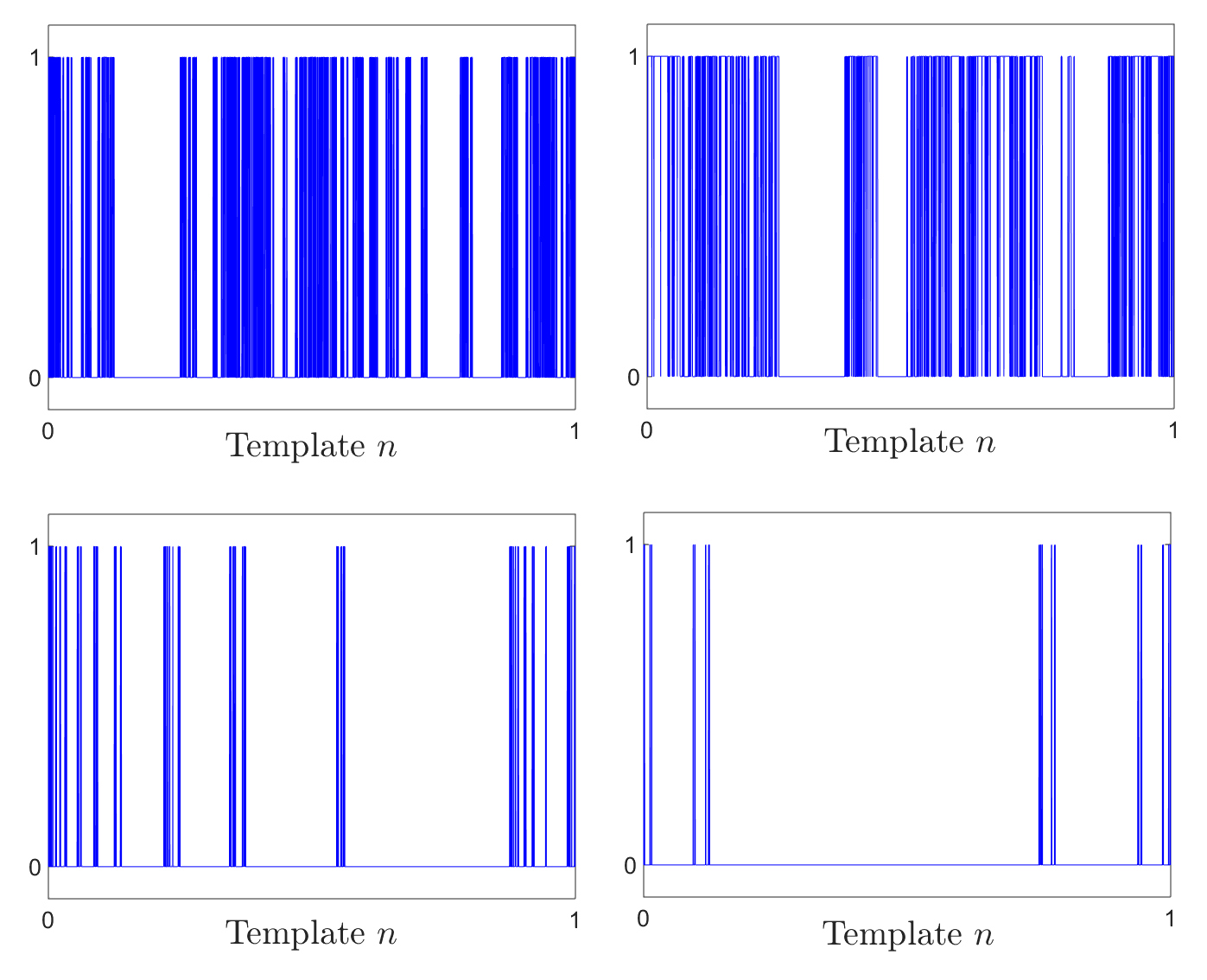}
\end{center}
\caption{\emph{{\bf Fixed map Mandelbrot sets}, for different parameter pairs $(c_0,c_1) \in \mathbb{C}$. {\bf A.} $c_1=-0.117-0.76i$, $c_0=-0.5622-0.62i$;  {\bf B.} $c_1=-0.5622-0.62i$, $c_0=-0.5622-0.62i$; {\bf C.} $c_1=-0.117-0.76$, $c_0=-0.75$; {\bf D.} $c_1=-0.75$, $c_0=-0.117-0.856i$. We created binary expansions of lentgh $L=15$ for the numbers in the unit interval $[0,1]$, which we then used as the symbolic templates for the iteration process.}}
\label{map_mand}
\end{figure}

In this section, we suggest some alternative ways of viewing the fixed map Mandelbrot set, and of observing how the template affects the properties of the orbit $o_{\bf s}(0)$, for a fixed parameter pair $(c_0,c_1) \in \mathbb{C}$. As a book keeping method for truncated templates, we propose to consider each binary sequence of length $L$ as equivalent to the binary representation of a real number $0 \leq n \leq 1$ with precision $2^{-L}$. In other words, we can consider, for each $0 \leq a \leq 1$, its binary expansion up to its $L$-th binary digit, always choosing the infinite expansion (when it is the case). For fixed $c_0$ and $c_1$, we can use each expansion as the symbolic template, and check whether the orbit $o_{\bf s}(0)$ is bounded or not. We can construct the function: $F \colon [0,1] \rightarrow \{ 0,1 \}$, given by $F(n)=1$, if $o_{\bf s}$ is bounded, and $F(n)=0$, if $o_{\bf s}(0)$ is not bounded. In Figure~\ref{map_mand} we show a few such representations, for various choices of the complex parameter pair $(c_0,c_1)$, connecting the discontinuous points in order to make the structure of the set more visible. It is clear that the behavior of the fixed map Mandelbrot set $F^{-1}(1)$ depends crucially on the parameter values. However, an analysis of its topological properties (accumulation points, density etc), or its measure in $[0,1]$ remains open for further studies.

\begin{figure}[h!]
\begin{center}
\includegraphics[width=0.8\textwidth]{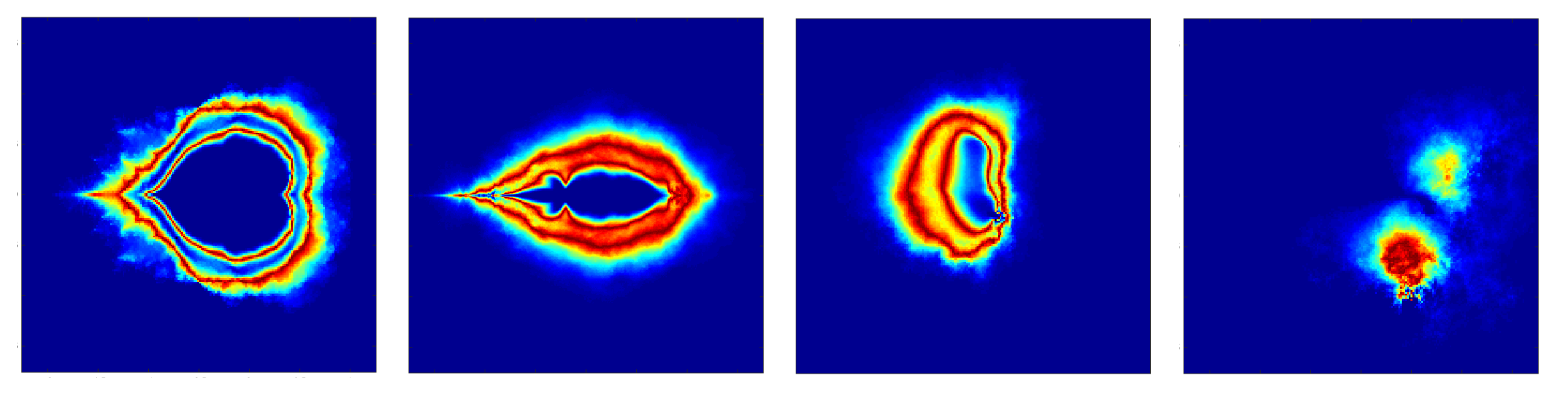}
\end{center}
\caption{\emph{{\bf Hybrid Mandelbrot sets} for different values of $c_0$: {\bf A.} $c_0=0$; {\bf B.} $c_0=-0.75$; {\bf C.} $c_0=0.375+0.333i$; {\bf D.} $c_0=i$. For each $c_1$ (represented in the complex plane), we used colors to illustrate how many templates of length $L=20$ lead to a bounded orbit $o_{\bf s}(0)$. The color spectrum goes from blue (low) to red (high).}}
\label{hybrid}
\end{figure}

Another interesting representation can be obtained as a hybrid of fixing the map and fixing the template. With a fixed $c_0$, we can measure, for each different value of $c_1 \in \mathbb{C}$, how likely it is (i.e., for how many of all templates of a certain length) that the orbit $o_{\bf s}(0)$ is bounded. In Figure~\ref{hybrid} we show a few such slices in the $c_1$ complex plane, with the colors representing the number of templates for which $o_{\bf s}(0)$ is bounded.

\section{Discussion}
\label{discussion}

\subsection{Comments and future work}

In this paper, we described a problem which is, to the best of our knowledge, a new extension of a well known traditional problem in discrete dynamics. We defined template iterations of two quadratic maps, and we used standard numerical algorithms to investigate the Julia and Mandelbrot sets of template systems. While this study is a first step in establishing the framework and phrasing some conjectures, a lot of theoretical work remains to be done, possibly exploring a brand new set of questions in discrete dynamics, which are likely to require new, different methods than the ones used in the traditional context of Fatou an Julia. 

The dynamic behavior in the case of one iterated map, and subsequently the connectivity type of classical Julia sets have been shown to be completely determined by the map's critical orbits. While related approaches may still be appropriate to use in the case of alternating maps~\cite{danca2009alternated}, and even in the case of higher template periods, this direction is not helpful when attempting to study the case of nonperiodic templates (which is in fact more plausible as a model for applications in the natural sciences). That is because the concepts of periodic or critical orbits are in this context no longer well defined, hence one has to start building a theoretical framework from its bare bones. In our current work we are exploring new methods, in particular expanding the question in a more abstract framework, obtained by interpolating continuously in the space of binary templates.

In the meantime, in our computational work we are continuing to investigate extensions of the Julia set, in particular in the context of networks of interconnected complex quadratic maps. While pairs of coupled logistic maps have been studied before in both real and complex case,   we are interested to study the coupled behavior in higher dimensional networks, with possible applications to understanding dynamics in neural networks. We are focused in particular on understanding the effects of the network architecture (e.g., graph Laplacian) and the dynamic properties of the ensemble as a whole, and on the topology of the ``Julia set'' of the networked system.

\subsection{Applications to genetics}

When a cell divides, it has to copy and transmit the exact same sequence of billion nucleotides to its daughter cells. While most DNA is typically copied with high fidelity (polymerase enzymes are amazingly precise when performing DNA synthesis), errors are a natural part of DNA replication, with rates of about 1 per $10^5$ (polymerases sometimes inserting too many or too few, or erroneous nucleotides into a sequence).  Human diploid cells have 6 billion base pairs, and each cell division makes about 120,000 errors~\cite{pray2008dna}.

Cells have evolved highly sophisticated DNA repair processes, aimed to promptly fix most of these errors. Some errors are corrected right away, during replication, through a repair process known as proofreading. Proofreading fixes about 99\% of the errors, but that's still not good enough for normal cell functioning. Some errors are corrected after replication, in a process called mismatch repair. Incorrectly paired nucleotides that still remain following mismatch repair become permanent mutations after the next cell division: once established, the cell no longer recognizes them as errors~\cite{pray2008dna}, passing them on to next generations of cells and (if the errors occur in gemetes) even to next generations of the organism.

When the genes for the DNA repair enzymes themselves become mutated (the iterated function changes in the long-term), mistakes begin accumulating at a much higher rate. Mutation rates vary substantially among taxa, and even among different parts of the genome in a single organism. Scientists have reported mutation rates as low as 1 per $10^{6}$-$10^9$ nucleotides, mostly in bacteria, and as high as 1 per $10^2$-$10^3$ nucleotides in humans~\cite{johnson2000fidelity}. Cells accumulate mutations as they divide. Even mutation rates as low as $10^{-10}$ can accumulate quickly over time, particularly in rapidly reproducing organisms like bacteria.

Our system can be viewed as a theoretical framework for studying iterated replication mechanisms which are subject to errors at each iteration step.  In genetics, polymerases replicate identically the DNA strand at division, which in turn governs the development of the cell, and is passed on at the next cell division. In our model, the original system/cell (in our case, the complex $z$-plane) is programmed to evolve according to a certain sequence of steps, leading to emergence of some features and extinction of others. For example, an initial $\xi_0$ which iterates to $\infty$ may represent a cell feature which becomes unsustainable after a number of divisions, while an initial $\xi_0$ which is attracted to a simple periodic orbit may represent a feature which is too simple to be relevant or efficient for the cell. Then the points on the boundary between these two behaviors (i.e., the Julia set) may be viewed as the optimal features, allowing the cell to perform its complex function. An error at the level of the iteration function at one particular iteration step is equivalent to a mutation that occurred at one of the cell division steps. The new cell/complex plane is then used as template for the next iteration/division; one can study how the features of the cell are affected in the long term, when such an error passes undetected by the repair mechanisms. In our paper, we considered situations where such an occurrence is singular, random/occasional or periodic. It is clear that mutations accumulated over a long period of time may lead to serious changes in the structure of the later cells (different topological properties of the Julia set). Our model also addresses the timing when the errors occur, and illustrates how a mutation in the early iterative stages can lead to substantially more dramatic consequences on the result (Julia set) than the same error if it happens later in the process. 

The construction of mathematical models to help understand DNA replication and repair would be highly desired, since these are crucially important and complex mechanisms to study and understand. In eukariotic cells, accumulating mutations can lead to cancer. However, if DNA replication were perfect (mutation-free), there would be no genetic variation. Therefore, successful organisms had to construct optimal mechanisms, providing efficient DNA repair, but also enough variability for evolution to continue. A mathematical framework would be ideal for posing and contextualizing such questions.\\


\clearpage

\thispagestyle{empty}

\begin{landscape}

\setlength{\textheight}{10in} \setlength{\voffset}{-0.8in}

\section*{Appendix 1: Effect of a propagating error along the binary template }

\vspace{7mm}
\begin{figure}[h!]
\begin{center}
\includegraphics[width=1.5\textwidth]{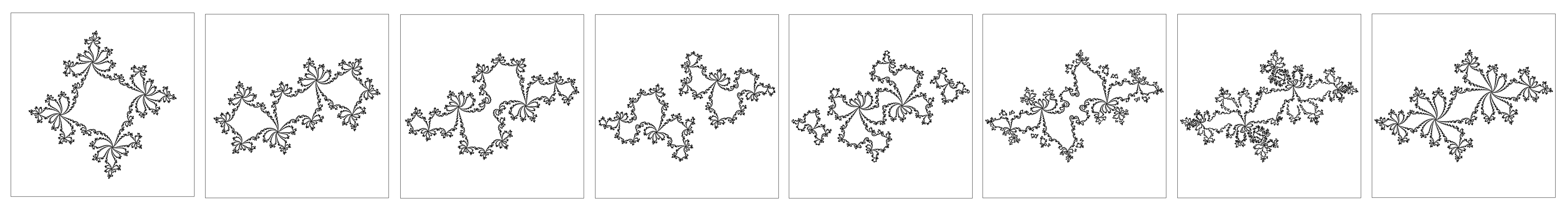}
\end{center}
\caption{\emph{{\bf Effect of error propagation on the Julia set} for the desired function with parameter $c_1=-0.62-0.432i$ (whose classical Julia set is shown in Figure~\ref{Julia_classic}c) and error parameter $c_0=0$ (whose classical Julia set is the unit circle). The perturbation of $f_{c_1}$ to $f_{c_0}$ was introduced successively at the iterations $k=1, 2, 3, 4, 5, 10, 30$ and $200$ in a truncated template of length $N=200$ (each Julia set is represented in one of the figure panels, from left to right).}}
\label{error_anim_c3}
\end{figure}

\vspace{1cm}
\begin{figure}[h!]
\begin{center}
\includegraphics[width=1.5\textwidth]{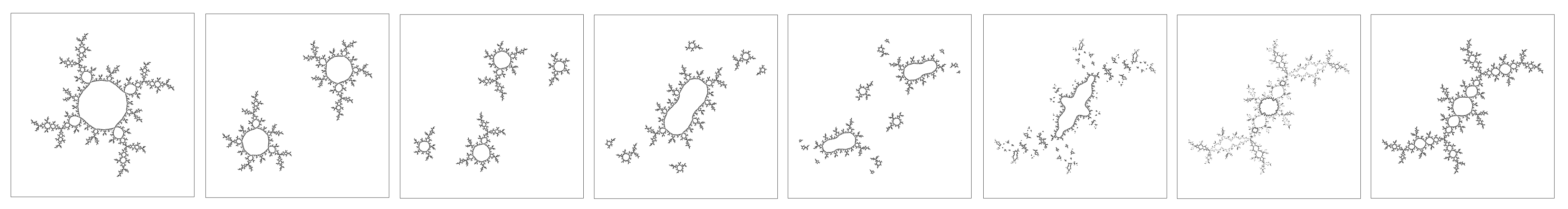}
\end{center}
\caption{\emph{{\bf Effect of error propagation on the Julia set} for the desired function with parameter $c_1=-0.117-0.856i$ (whose classical Julia set is shown in Figure~\ref{Julia_classic}e) and error parameter $c_0=0$ (whose classical Julia set is the unit circle). The perturbation of $f_{c_1}$ to $f_{c_0}$ was introduced successively at the iterations $k=1, 2, 3, 4, 5, 10, 30$ and $200$ in a truncated template of length $N=200$ (each Julia set is represented in one of the figure panels, from left to right).}}
\label{error_anim_c7}
\end{figure}

\end{landscape}

\clearpage

\thispagestyle{empty}

\begin{landscape}

\setlength{\textheight}{10in} \setlength{\voffset}{-0.8in}

\begin{figure}[h!]
\begin{center}
\includegraphics[width=1.5\textwidth]{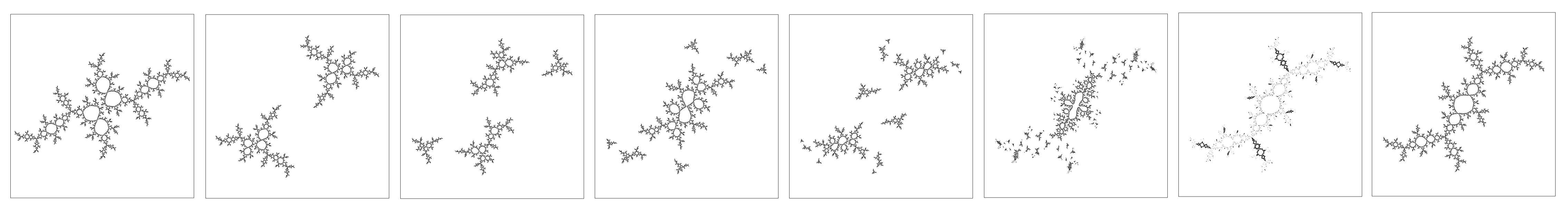}
\end{center}
\caption{\emph{{\bf Effect of error propagation on the Julia set} for the desired function with parameter $c_1=-0.117-0.856i$ (whose classical Julia set is shown in Figure~\ref{Julia_classic}d) and error parameter $c_0=-0.5622-0.62i$ (whose classical Julia set is shown in Figure~\ref{Julia_classic}e). The perturbation of $f_{c_1}$ to $f_{c_0}$ was introduced successively at the iterations $k=1, 2, 3, 4, 5, 10, 30$ and $200$ in a truncated template of length $N=200$ (each Julia set is represented in one of the figure panels, from left to right).}}
\label{error_anim_c7_c4}
\end{figure}

\vspace{1cm}
\begin{figure}[h!]
\begin{center}
\includegraphics[width=1.5\textwidth]{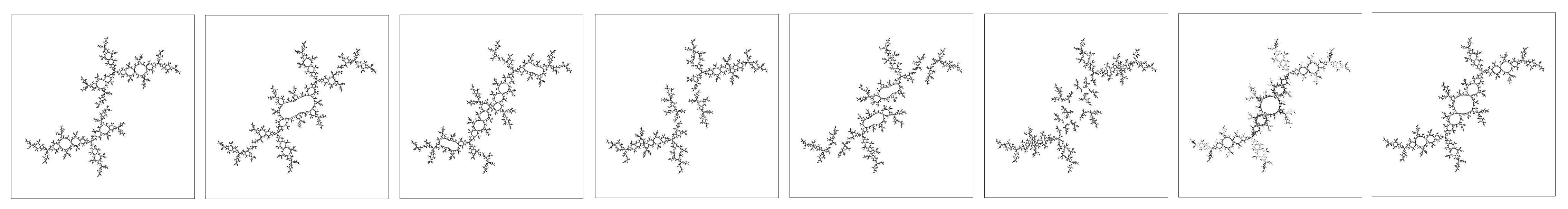}
\end{center}
\caption{\emph{{\bf Effect of error propagation on the Julia set} for the desired function with parameter $c_1=-0.117-0.856i$ (whose classical Julia set is shown in Figure~\ref{Julia_classic}d) and error parameter $c_0=c_1 + \varepsilon$, where $\varepsilon = 0.1+0.1i$ (a small complex perturbation of $c_1$). The perturbation of $f_{c_1}$ to $f_{c_0}$ was introduced successively at the iterations $k=1, 2, 3, 4, 5, 10, 30$ and $200$ in a truncated template of length $N=200$ (each Julia set is represented in one of the figure panels, from left to right).}}
\label{error_anim_c7_perturb}
\end{figure}

\end{landscape}

\clearpage

\thispagestyle{empty}

\begin{landscape}

\setlength{\textheight}{10in} \setlength{\voffset}{-0.8in}

\section*{Appendix 2: Template Mandelbrot slices}

\begin{figure}[h!]
\begin{center}
\includegraphics[width=1.1\textwidth]{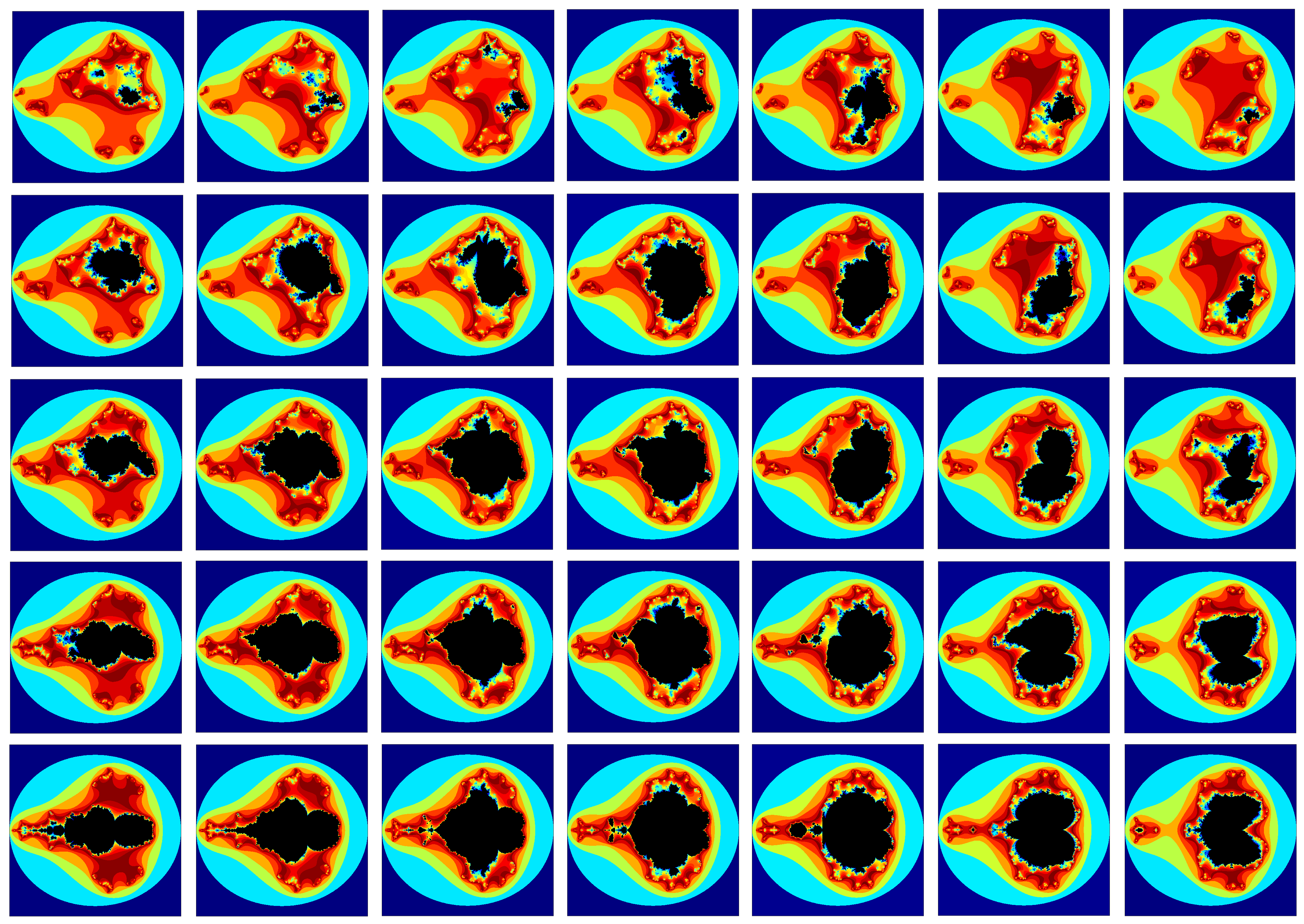}
\end{center}
\caption{\emph{ {\bf Mandelbrot slices for a lattice of $c_0$ or step 0.2}, ranging from -0.6 to 0.6 in the real direction, and from 0 to 0.8 in the imaginary direction. The slices are symmetric with respect to the real axis (not snown).}}
\label{Mand_all_011}
\end{figure}

\end{landscape}

\end{document}